\documentclass[12pt,reqno]{article}
\usepackage{amsmath,amsthm,amsfonts,amssymb,amscd}
\usepackage[dvips]{graphicx}
\usepackage{psfrag}

\input{psfig}
\theoremstyle{plain}
\newtheorem{thm}{Theorem}[section]
\newtheorem{rk}[thm]{Remark}
\newtheorem{prop}[thm]{Proposition}

\newtheorem{lemma}[thm]{Lemma}

\newtheorem{defi}[thm]{Definition}

\newtheorem{maintheorem}{Theorem}

\newcommand{\interior}{\operatorname{int}}

\newcommand{\si}{\operatorname{Sing}}

\newcommand{\vol}{\operatorname{vol}}

\newcommand{\dom}{\operatorname{Dom}}

\newcommand{\diam}{\operatorname{diam}}

\title{Topological dimension of
singular-hyperbolic attractors}

\author{C. A. Morales\thanks{2000 MSC: Primary 37D30,
Secondary 37C45.
{\em Key words and phrases}:
Attractor, Partially Hyperbolic,
Topological Dimension.
Partially supported by CNPq, FAPERJ and PRONEX-Dyn. Sys./Brazil.}}

\begin{document}

\maketitle

\begin{abstract}
An {\em attractor} is a transitive set
of a flow to which all positive orbit
close to it converges.
An attractor is {\em singular-hyperbolic}
if it has singularities (all hyperbolic)
and is partially hyperbolic with volume expanding
central direction \cite{MPP}.
The geometric Lorenz attractor
\cite{GW} is an example of a singular-hyperbolic attractor with topological dimension $\geq 2$.
We shall prove that {\em all} singular-hyperbolic attractors
on compact $3$-manifolds
have topological dimension $\geq 2$.
The proof uses the methods in \cite{MP}.
\end{abstract}

\section{Introduction}
\label{s-intro}

This paper is concerned with
the topological dimension of attractors for flows on
compact manifolds.
By {\em attractor} we mean a transitive set
of the flow to which all positive orbit
close to it converges.
The attractors under consideration
will be {\em singular-hyperbolic} in the sence that
they have singularities (all hyperbolic)
and
are partially hyperbolic with
volume expanding central direction \cite{MPP}.
In particular, the singular-hyperbolic attractors
are volume hyperbolic sets as defined in \cite{B}.
The geometric Lorenz attractor
is an example of a singular-hyperbolic
attractor
with topological dimension $\geq 2$.
We shall prove that {\em all}
singular-hyperbolic attractors
on compact $3$-manifolds have
topological dimension $\geq 2$.
The proof uses the methods developed in \cite{MP}.
Let us state our result in a precise way.

Hereafter $X$ will be a $C^1$ vector field
on a compact manifold $M$. The flow of $X$
is denoted by $X_t$, $t\in I\!\! R$.
Given $p\in M$ we define
$\omega(p)=\omega_X(p)$, the $\omega$-limit set of $p$,
as the accumulation point set of the positive orbit
of $p$.
The {\em $\alpha$-limit set} of $p$ is
the set $\alpha(p)=\alpha_X(p)=\omega_{-X}(p)$.
A compact invariant set $\Lambda$ of $X$ is {\em transitive}
or {\em attracting} depending on whether
$\Lambda=\omega(p)$ or $\cap_{t>0}X_t(U)$
for some $p\in \Lambda$ or
some compact neighborhood $U$ of $\Lambda$ respectively.
An {\em attractor} is a
transitive attracting sets.
A closed orbit of $X$ is either
periodic or singular.
A singularity of $X$ is
{\em hyperbolic} if
none of its eigenvalues
have zero real part.

A compact invariant set
$\Lambda$ of $X$ is {\em partially hyperbolic} \cite{HPS}
if there are an invariant splitting
$T\Lambda=E^s\oplus E^c$
and positive constants $K,\lambda$
such that:

\begin{enumerate}
\item
{\em $E^s$ is contracting}, namely
$$
\mid\mid DX_t/E^s_x\mid \mid\leq
K e^{-\lambda t},
\,\,\,\,\,\,\forall x\in \Lambda,\,\,\forall
t>0.
$$
\item
{\em $E^s$ dominates $E^c$}, namely
$$
\mid \mid DX_t/E^s_x\mid\mid
\cdot
\mid\mid
DX_{-t}/E^c_{X_{t}(x)}
\mid\mid
\leq K e^{-\lambda t},
\,\,\,\,\,\,\forall x\in \Lambda,\,\,\forall
t>0.
$$
\end{enumerate}

The central direction $E^c$ of
$\Lambda$ is said to be {\em volume expanding}
if the additional condition
$$
\mid Det(DX_t/E^c_x)\mid
\leq K e^{-\lambda t}
$$
holds $\forall x\in \Lambda$, $\forall t>0$
where
$Det(\cdot)$ means the jacobian.
The above splitting $E^s\oplus E^c$ will be refered to
as a {\em $(K,\lambda )$-splitting}
in the Appendix.

\begin{defi}(\cite{MPP})
\label{d1}
An attractor is {\em singular-hyperbolic}
if it has singularities (all hyperbolic) and is
partially hyperbolic with volume expanding central direction.
\end{defi}

\begin{defi}(\cite{HW})
\label{top-dim}
The {\em topological dimension} of a
space $E$ is either $-1$
(if $E=\emptyset$) or
the last integer
$k$ for which every point has arbitrarily small
neighborhoods whose boundaries
have dimension less than $k$.
\end{defi}

The relation between dynamics and topological
dimension was considered
for hyperbolic systems \cite{H,P,B1,B2}; for expansive systems \cite{KS,M};
and for singular-hyperbolic systems \cite{Mo}.
The result below generalizes to singular-hyperbolic attractors a well known property of both hyperbolic strange attractors and geometric Lorenz attractors.

\begin{maintheorem}
\label{thA}
Singular-hyperbolic attractors on compact $3$-manifolds
have topological dimension $\geq 2$.
\end{maintheorem}

The idea of the proof is the following.
Let $\Lambda$ be a singular-hyperbolic
attractor of a flow $X$ on a compact $3$-manifold $M$.
It follows from \cite{MPP} that
all the singularities $\sigma\in \Lambda$ are
{\em Lorenz-like},
namely the eigenvalues
$\lambda_1,\lambda_2,\lambda_3$ of $\sigma$
are real and satisfy
$\lambda_2<\lambda_3<0<-\lambda_3<\lambda_1$.
The flow nearby $\sigma$ can be described
using the Grobman-Hartman Theorem \cite{dMP}.
In particular, a Lorenz-like singularity
exhibits two 
{\em singular cross-sections} $S^t,S^b$
and two {\em singular curves} $l^t,l^b$
(\cite{MP}).
A {\em singular cross section} of $\Lambda$
is by definition a disjoint collection of singular cross sections
$S^t,S^b$ (as $\sigma$ runs over all the singularities
of $\Lambda$) whose horizontal boundaries
does not intersect $\Lambda$.
The {\em singular curve} of $S$ is the
union $l$ of the respective singular curves $l^t,l^b$.
A {\em singular partition} of $\Lambda$
will be a compact neighborhood $O$ of
$\Lambda\cap l$ in $S$, for some
singular cross section $S$ of $\Lambda$,
such that $\Lambda\cap l$ does not intersect
the boundary of $O$ and
every regular orbit of $\Lambda$ intersect $O$.
The {\em size} of the singular partition $O$
is the minimal $\epsilon>0$ such that there is
an invariant cone field
in $O$ (for the return map
$\Pi:\dom(\Pi)\subset O\to O$)
on which the derivative of $\Pi$ has
expansion rate bigger than $\epsilon^{-1}$.
In Proposition \ref{p1} we prove that
one-dimensional singular-hyperbolic attractors on compact $3$-manifolds have singular partition
with arbitrarily small size.
The proof of this proposition uses
the Lemmas 7.5 and 7.6 in \cite{MP}. These lemmas
will be proved in the Appendix for the sake of completeness.
In Theorem \ref{p2} we shall prove that
singular-hyperbolic attractors $\Lambda$ on compact $3$-manifolds
cannot have singular partitions with arbitrarily small size.
Theorem \ref{thA} will follow from
Proposition \ref{p1} and Theorem \ref{p2}.

\section{Proof}

We start with some definitions.
Hereafter $\Lambda$ is a singular-hyperbolic attractor
of a $C^1$ flow $X$ on a compact $3$-manifold $M$.
Since Lorenz-like singularities $\sigma$ are hyperbolic
they are equipped with three invariant manifolds
$W^s_X(\sigma),W^u_X(\sigma),W^{ss}_X(\sigma)$
each one tangent at $\sigma$ to the eigenspace corresponding
to $\{\lambda_2,\lambda_3\},\{\lambda_1\},\{\lambda_2\}$
respectively.
It follows from
\cite{MPP} that every singularity $\sigma$ of $X$ in $\Lambda$ is Lorenz-like
and satisfies
$
\Lambda\cap W^{ss}_X(\sigma)=\{\sigma\}
$.
The classical Grobman-Hartman Theorem \cite{dMP}
gives the description of the flow nearby $\sigma$.
This is done at Figure \ref{f.1}.
Note that $W^{ss}_X(\sigma)$ separates $W^s_X(\sigma)$
in two connected components denoted the top and the bottom respectively.
In one of these components, say the top one, we consider
a cross-section $S^t=S^t_\sigma$ together with a curve $l^t=l^t_\sigma$ as in Figure \ref{f.1}.
Similarly we consider a cross-section $S^b=
S^b_\sigma$ and a curve $l^b=l^b_\sigma$ located
in the bottom component of $W^s_X(\sigma)$. See Figure \ref{f.1}.
Both $S^*$ (for $*=1,2$) are homeomorphic to $[0,1]\times [0,1]$.
$S^*$ can be chosen in a way that $l^*$ is contained 
in $W^s_X(\sigma)\setminus W^{ss}_X(\sigma)$.
The positive flow lines of $X$ starting at $S^t\cup S^b\setminus(l^t\cup l^b)$
exit a small neighborhood of $\sigma$ passing through the cusp region
as indicated in Figure \ref{f.1}.
The positive orbits starting at $l^t\cup l^b$ goes directly to $\sigma$.
We note that the boundary of $S^*$ is formed by four curves,
two of them transverse to $l^*$ and two of them
parallel to $l^*$.
The union of the
curves in the boundary of $S^*$
which are parallel (resp. transverse)
to $l^*$ is denoted by $\partial^v S^*$ (resp.
$\partial^h S^*$).
The interior (as a submanifold) of
$S^*$ is denoted by $Int(S^*)$.

\begin{rk}
\label{rk}
An immediate consequence
of $\Lambda\cap W^{ss}_X(\sigma)=\{\sigma\}$
is the following.
Let $\sigma$ be a singularity of $X$ in $\Lambda$.
Then there are cross-sections $S^t,S^b$ as above
arbitrarily close to
$\sigma$ such that
$\Lambda\cap \partial^h S^*=\emptyset$
($*=t,b$).
Since the two boundary points
of $l^*$ are in $\partial^h S^*$
we have that $\Lambda\cap l^*\subset Int(S^*)$.
\end{rk}

\begin{defi}
\label{AL}
We shall call the cross sections $S^t,S^b$
as {\em singular cross sections} associated to $\sigma$.
The curves $l^t,l^b$ are called
{\em singular curves} of $S^t,S^b$ respectively.
A {\em singular cross section} of $\Lambda$ is
a finite disjoint collection
$
\{S^t_\sigma,S^b_\sigma:
\sigma$ is a singularity of $X$ in $\Lambda\}$
satisfying $\Lambda\cap \partial^h S=\emptyset$.
The {\em singular curve} of $S$ is the
associated collection of singular curves
$
l=
\{l^t_\sigma,l^b_\sigma:
\sigma$ is a singularity of $X$ in $\Lambda\}$.
\end{defi}

%\begin{figure}[htv] 
%\centerline{
%\psfig{figure=ly=at1.eps,height=3in}}
%\caption{\label{f.1} Lorenz-like singularity.}
%\end{figure}

Hereafter we denote by $T_\Lambda M=E^s_\Lambda\oplus E^c_\Lambda$ the singular-hyperbolic
splitting of $\Lambda$.
The contracting direction $E^s$ is
one-dimensional and contracting.
So,
$E^s_\Lambda$ can be extended
to an invariant contracting splitting 
$E^s_{U(\Lambda)}$ on a neighborhood $U(\Lambda)$ of $\Lambda$. The standard Invariant Manifold Theory
\cite{HPS} implies that
$E^s_{U(\Lambda)}$
is tangent to a continuous
foliation ${\cal F}$ on $U(\Lambda)$. 
If $S$ is a singular cross-section
contained in $U(\Lambda)$,
we denote by ${\cal F}^{S}$ the foliation of $S$ obtained  projecting 
${\cal F}$ into $S$ along $X$.
The space of leaves of ${\cal F}^S$ will be denoted by
$I^S$.
We extend $E^c_\Lambda$ continuously to a subbundle $E^c_{U(\Lambda)}$ 
of $T_{U(\Lambda)}M$.
{\em In what follows we fix such a neighborhood $U(\Lambda)$ of $\Lambda$}.

\begin{rk}
\label{rk'}
It is possible to choose $S$ arbitrarily close to the singularities of $\Lambda$ in a way that
$l$ is a finite union of leaves of ${\cal F}^S$ and
$I^S$ is a finite disjoint union of compact intervals.
\end{rk}

The following lemma is a direct consequence of standard
argument involving topological dimension.
We prove it here for the sake of completeness.

\begin{lemma}
\label{hw1}
Let $S$ a singular cross-section
and $l$ be its associated singular curve.
If $\Lambda$ is one-dimensional,
then there is a compact neighborhood
$O$ of $\Lambda\cap l$ in $S$
whose boundary $\partial O$ satisfies
$\Lambda\cap \partial O=\emptyset$.
\end{lemma}

\begin{proof}
Note that
$\Lambda\cap \partial^hS=\emptyset$ since
$S$ is a singular cross-section. As noted in Remark \ref{rk} one has
$\Lambda\cap l\subset Int(S)$.
Fix $x\in \Lambda\cap l$.
Then $x\in Int(S)$.
Because $\Lambda$ is one dimensional we have
that $\Lambda\cap S$ is zero dimensional \cite{HW}.
Then, by the definition of the topological dimension,
one can find an open set
$S_x$ of $\Lambda\cap S$ containing $x$
such that $\partial S_x=\emptyset$.
Note that the topology in
$\Lambda\cap S$ is the one induced by
$S$.
In follows that
$S_x=(\Lambda\cap S)\cap O_x$
for some open set $O_x$ of $M$.
Since $S$ is transversal to $X$
we can choose $O_x$ such that
$\partial S_x=(\Lambda\cap S)\cap \partial O_x$
(for this we can use the Tubular Flow-Box Theorem
\cite{dMP}).
It follows that
$$
(\Lambda\cap S)\cap \partial O_x=\emptyset.
$$
On the other hand,
$\Lambda\cap l$ is compact in $S$ and
$\{S\cap O_x:x\in \Lambda\cap l\}$ is
an open covering of $\Lambda\cap l$.
It follows that there is a finite
subcollection of $\{S\cap O_x:x\in \Lambda\cap l\}$
covering $\Lambda\cap l$.
Denote by $O$ the union of the closures (in $S$)
of the elements
of such a subcollection.
It follows that $O$ is a compact neighborhood
of $\Lambda\cap l$ in $S$.
Since $O$ is a finite union of
$S\cap O_x$'s satisfying
$(\Lambda\cap S)\cap \partial O_x=\emptyset$ we have that $\Lambda\cap \partial O=\emptyset$.
This proves the lemma.
\end{proof}

Herefter $O$ is a set
contained in a singular cross section $S$.
Clearly $O$ defines a return map
$$
\Pi:\dom (\Pi)\subset O\to O
$$
given by
$$
\Pi(x)=X_{t(x)}(x),
$$
where $\dom(\cdot)$ denotes the domain and $t(\cdot)$ denotes the return time.

\begin{rk}
\label{rk''}
Note that $\Pi$
may be discontinuous in $\Pi^{-1}(\partial O)$.
However if $x\in \Pi^{-1}(Int(O))$
then $\Pi$ is $C^1$ in an open neighborhood
of $x$ contained in $Int(O)$.
This is an immediate consequence of the
Tubular Flow-Box Theorem.
\end{rk}

We denote by $TO$
the tangent space of $O$ relative $S$.
If $x\in M$ we denote by
$\angle(v_x,w_x)$ the
tangent of the
angle between $v_x,w_x\in T_xM$.
If $L_x$ is a linear subspace of $T_xM$, we define
$$
\angle(v_x,L_x)=\inf_{w_x\in L_x}\angle(v_x,w_x).
$$
Given $\alpha>0$ we define the {\em cone}
$$
C_\alpha(L_x)=\{v_x\in T_xM:\angle(v_x,L_x)\leq \alpha\}.
$$
If $L:x\in Dom(L)\to L_x$ is a map and $\alpha>0$ we define
the {\em cone field}
$$
{\cal C}_\alpha(L)=
\{C_\alpha(L_x):x\in Dom(L)\}.
$$
The case $L=E^c$ will be interesting.
The definition below is a minor modification
of the corresponding definition in
\cite{MP}.
If $x\in M$ we denote
$X_{I\!\! R}(x)$ the full orbit of $x$.

\begin{defi}
\label{sp}
A {\em singular partition} of $\Lambda$ is
a set $O$ satisfying the following properties:

\begin{enumerate} 
\item
There is a singular cross-section $S$ such that
$O\subset Int(S)$ is a compact neighborhood
of $\Lambda\cap l$. 
\item
$\Lambda\cap \partial O=\emptyset$.
\item
$
\si_X(\Lambda)=\{q\in \Lambda:
X_{I\!\! R}(q)\cap O=\emptyset\}.
$
\end{enumerate}

{\em The size of $O$} is the minimal number
$\epsilon>0$ for which there is $\alpha>0$ such that
the cone field
${\cal C}_\alpha(E^c)$
satisfies :

\begin{description}
\item{4.}
If $x\in \dom(\Pi)$, then
$$
D\Pi(x)\left(
C_\alpha(E^c_x)\cap TO_x
\right)\subset 
\interior\left(
C_{\frac{\alpha}{2}}(E^c_{\Pi(x)})\cap TO_{\Pi(x)}
\right).
$$
\item{5.}
If $x\in \dom(\Pi)$ and $v_x\in C_\alpha(E^c_x)\cap TO_x$, then
$$
\| D\Pi(x)(v_x)\| \geq \epsilon^{-1}\| v_x \|.
$$
\item{6.}
$
\inf\{\angle(v_x,E^s_x): x\in O, v_x\in
C_\alpha(E^c_x)\cap TO_x\}>0
$.
\end{description}
\end{defi}

The following proposition
studies the existence of singular partition
with arbitrarily small size
for certain one-dimensional singular-hyperbolic sets.
Its proof uses the
methods developed in \cite{MP}.
We let
$\si_X(\Lambda)=\{\sigma_1,\cdots ,\sigma_k\}$ be
the set of singularities of $X$ in $\Lambda$.

\begin{prop}
\label{p1}
One-dimensional singular-hyperbolic
attractors on compact
$3$-manifolds have singular
partitions with arbitrarily small size.
\end{prop}

\begin{proof}
Let $\Lambda$ be a singular-hyperbolic attractor
of a $C^1$ flow $X$ on a compact $3$-manifold $M$.
We shall assume that $\Lambda$ has topological
dimension
$1$. We shall prove that
$\Lambda$ has singular partition with arbitrarily
small size $\epsilon>0$.
For this we proceed as follows.
Since $\Lambda$ has topological dimension $1$
we have that $\Lambda$ cannot contain
hyperbolic sets (the unstable manifold
of a hyperbolic set in
$\Lambda$ would be two-dimensional and contained in $\Lambda$).
It follows that $\omega(x)$ cannot be hyperbolic
for all $x\in \Lambda$. By \cite{MPP'}
if ${\cal L}=\omega(x),\alpha(x)$ then
\begin{equation}
\label{**}
{\cal L}\cap Sing_X(\Lambda)\neq\emptyset,
\,\,\,\,\,\forall x\in \Lambda.
\end{equation}

Choose $\alpha>0$ such that
$$
\inf
\{\angle(v_x,E^s_x):x\in U(\Lambda),v_x\in C_\alpha(E^c_x)\}>0.
$$
By 
\cite[Lemma 7.5]{MP} (see Lemma \ref{AP1})
we can find a neighborhood
$U_\alpha\subset U(\Lambda)$ of
$\Lambda$ and positive constants
$T_\alpha,K_\alpha,\lambda_\alpha$
such that the following properties hold:

\begin{description}
\item{(P1).}
If $x\in U_\alpha$ and $t\geq T_\alpha$,
then
$$
DX_t(x)(C_\alpha(E^c_x))
\subset C_{\alpha/2}(E^c_{X_t(x)}).
$$
\item{(P2).}
If $x\in U_\alpha$ is regular,
$X(x)\in C_\alpha(E^c_x)$, $t\geq T_\alpha$ and $v_x\in
C_\alpha(E^c_x)$ is orthogonal to
$X(x)$, then
$$
\mid\mid
P^t_x(v_x)
\mid\mid\cdot
\mid\mid
X(X_t(x))
\mid\mid
\geq
K_\alpha e^{\lambda_\alpha t}
\cdot
\mid\mid v_x
\mid\mid
\cdot
\mid\mid
X(x)
\mid\mid,
$$
\end{description}
where
$P^t_x$ denotes the {\em Poincar\'e flow
associated to $X$} (see \cite{D,MP} or the Appendix).

Once we fix $\alpha$ and $U_\alpha$ we apply
\cite[Lemma 7.6]{MP} (see Lemma \ref{AP2})
to find,
for every $i\in \{1,\cdots ,k\}$, a pair
of singular cross-sections
$S^{*,0}_i$ associated to $\sigma_i$
($*=t,b$) such that
$$
X(x)\in C_\alpha(x),\,\,\,\,\forall x
\in S^{*,0}_i.
$$
Define
$$
S^0=
\cup_{i=1}^k(S^{t,0}_i\cup S^{b,0}_i).
$$
It is clear that $S^0$ is a singular cross-section.
We denote by $l^0$ the singular curve of $S^0$.
Since $S^0$ is transversal to $X$
one can find a constant $D>0$
(depending on $S^0$) such that
$$
\frac{\mid\mid X(x)\mid\mid}{\mid\mid X(y)\mid\mid}
>D,
$$for all $x,y\in S^0$.
We choose $T_\epsilon>T_\alpha$
large enough so that
\begin{equation}
\label{*}
K_\alpha e^{\lambda_\alpha t} \cdot D>\epsilon^{-1}
\end{equation}
for all $t\geq T_\epsilon$.

For every $\delta>0$ we consider
a singular-cross section
$S^\delta\subset S^{0}$
($i=1,\cdots , k$ and $*=t,b$)
formed by small bands $S^{*,\delta}_i$
of diameter $2\delta$ around
the singular curve $l^{*}_i$ of
$S^{*,0}_i$. Note that the singular curve
of $S^\delta$ is $l^0$ (the one of $S^0$)
for all $\delta$.
Since $\Lambda$ is one-dimensional
Lemma \ref{hw1}
implies that $\forall x\in \Lambda\cap l^0$
there is a compact neighborhood
$O=O^\delta\subset S^{\delta}$
of $\Lambda\cap l^0$ such that
$\Lambda\cap \partial O=\emptyset$.
Note that $O$ is a singular partition of
$\Lambda$.
In fact, (1) and (2) of Definition \ref{sp} are
obvious. And
(3) of Definition \ref{sp}
follows from Eq. (\ref{**})
since $O$ is a compact neighborhood of
$\Lambda\cap l^0$.

Let us prove that if $\delta>0$ small enough, then
$O$ has size $\epsilon$.
For this we need to prove that
for $\delta$ small
the cone field ${\cal C}_\alpha(E^c)$
satisfies the properties (4)-(6) of Definition
\ref{sp}.
Let $\Pi:Dom(O)\subset O\to O$
be the return map induced by $X$ in $O$.
By definition
$\Pi(x)=X_{t(x)}(x)$
where
$t(x)$ is the return time of $x\in Dom(O)\subset O$
into $O$.
To calculate $D\Pi(x)$
we can assume without loss of generality
that
$S^0$ is orthogonal to $X$.
It follows that
$$
D\Pi(x)=P^{t(x)}_x
$$
for all $x\in Dom(\Pi)$.
Shrinking $\delta$
one has $t(x)>T_\epsilon$ for all
$x\in Dom(\Pi)\subset O$.
This allows us to apply the
properties (P1)-(P2) above.
In fact, since $
D\Pi(x)=P^{t(x)}_x
$ one has
$$
D\Pi(x)/C_\alpha(E^c_x)\cap TO_x
=
P^{t(x)}_x/C_\alpha(E^c_x)
\cap TO_x.
$$
Then,
Definition \ref{sp}-(4) follows from
(P1).
(P2) and Eq. (\ref{*})
imply
$$
\| D\Pi(x)(v_x)\|
=\|P^{t(x)}_x(v_x)\|
=
$$
$$
=\|
P^{t(x)}_x(v_x)\|\cdot
\| X(X_{t(x)}(x))\|
\cdot\| X(X_{t(x)}(x))\|^{-1}
\geq
$$
$$
\geq
K_\alpha e^{\lambda_\alpha t(x)}\cdot
\frac{\| X(x)\|}{\| X(X_{t(x)}(x))\|}\cdot
\| v_x\|
\geq K_\alpha e^{\lambda_\alpha t(x)}\cdot
D\cdot \| v_x\|
\geq \epsilon^{-1}\| v_x\|,
$$
$\forall x\in \dom(\Pi)$, $\forall v_x\in C_\alpha(E^c_x)
\cap TO_x$ because $X(x)\in C_\alpha(E^c_x)$
$\forall x\in S^0$ and $t(x)>T_\epsilon$
($\forall x\in \dom(\Pi)$).
This proves Definition \ref{sp}-(5).
Definition \ref{sp}-(6)
is a direct consequence of the choice
of $\alpha$.
The result follows.
\end{proof}

The following lemma will be used
to prove Theorem \ref{p2}.
Recall that if $S$ is a singular cross-section then
$S$ is endowed with a singular curve $l$.
If $O\subset S$ then $\Pi:\dom(\Pi)\subset O\to O$ denotes the return map associated to $O$.

\begin{lemma}
\label{util}
Let $\Lambda$ be a singular-hyperbolic attractor
of a flow $X$ on a compact $3$-manifold.
Let $O$ be a singular partition
of $\Lambda$. Then, there is
an open neighborhood $O'\subset O$
of $\Lambda\cap O$ such that:

\begin{enumerate}
\item
$O'\setminus l\subset \dom(\Pi)$.
\item
$\Pi$ is $C^1$ in $O'\setminus l$.
\item
$\Pi(O'\setminus l)\subset O'$.
\end{enumerate}
\end{lemma}

\begin{proof}
Because $\Lambda$ is an attractor we have that
the unstable manifold of any of its singularities
is contained in $\Lambda$.
In particular, every connected component
of $W^u_X(\sigma_i)\setminus \{\sigma_i\}$ is contained
in $\Lambda$ $\forall i$.
It follows from Definition \ref{sp}-(3) that
all such components intersect $O$.
By Definition \ref{sp}-(2)
such intersections can occur only in $Int(O)$.
This implies that there are small open
bands, centered at the singular curves
in $l$, whose union $V(l)$ satisfies
$V(l)\setminus l\subset \Pi^{-1}(Int(O))$.
As noted in Remark \ref{rk''} we have that
$\Pi$ is $C^1$ in $V(l)\setminus l$.
Again by
Definition \ref{sp}-(2)-(3)
one has $(\Lambda\cap O)\setminus V(l)\subset
\Pi^{-1}(Int(O))$.
So, by Remark \ref{rk''},
since $(\Lambda\cap O)\setminus V(l)$
is compact, there is an open neighborhhod
$V$ of $(\Lambda\cap O)\setminus V(l)$
contained in $\dom(\Pi)$
such that $\Pi$ is $C^1$ in $V$.
Observe that $V\cup V(l)$ is an open
neighborhood of $\Lambda\cap O$
such that $\Pi$ is $C^1$ in
$(V\cup V(l))\setminus l$.
On the other hand,
$\Lambda$ is an attractor by assumption.
Then, there is a neighborhood $U^*$
such that $X_t(U^*)\subset U^*$, $\forall t>0$.
Clearly one can choose $U^*$ to be
{\em arbitrarily close to $\Lambda$}.
In particular, 
$O':=O\cap U^*$ is contained in $V\cup V(l)$.
It follows that $O'\setminus l\subset \dom(\Pi)$
because $V\cup (V(l)\setminus l) \subset \dom(\Pi)$.
Because
$X_t(U^*)\subset U^*$ for all $t>0$
and the return time for the points in
$\dom(\Pi)$ is positive
we conclude that $\Pi(O'\setminus l)\subset O'$.
As $\Pi$ is $C^1$ in $(V\cup V(l))\setminus
l$ and $O'\subset V\cup V(l)$
we conclude that $\Pi$ is $C^1$
in $O'\setminus l$
This proves the result.
\end{proof}

\begin{thm}
\label{p2}
Singular-hyperbolic attractors
on compact $3$-manifolds cannot have
singular partitions with arbitrarily small size.
\end{thm}

\begin{proof}
Let $\Lambda$ be a singular-hyperbolic
attractor of a
$C^1$ flow $X$ on a compact $3$-manifold
$M$. By contradiction
we assume that $\Lambda$ has a singular partition
$O$ with arbitrarily small size $\epsilon>0$.
We fix $\epsilon\in (0,1/2)$.
We let $O'$ be the
open neighborhood obtained in Lemma \ref{util}
for $O$.
Hereafter we say that a $C^1$ connected
curve $c$ in $O$
is a {\em $C^u$-curve} if its tangent vector belongs
to the cone field $C_\alpha(E^c)$ at Definition \ref{sp}-(4).
Definition \ref{sp}-(4)
implies that $\Pi$ carries
$C^u$ curves in $\Pi^{-1}(Int(O))$
into $C^u$ curves in $O$ (see also Remark \ref{rk''}).
Definition \ref{sp}-(6) implies that
a $C^u$ curve in $O$ intersects
$l$ in at most one point $x_c$.
In that case $x_c$ divides
$c$ in two connected components
the largest one being denoted by $c^+$.
Clearly if $L(\cdot)$ denotes the lenght,
then
$$
L(c^+)\geq (1/2)L(c).
$$
Now, fix a $C^u$ curve $c_1\subset O'\setminus l$.
Define
$R=(2\epsilon)^{-1}$.
The choice of $\epsilon$
implies $\epsilon^{-1}>R>1$.
Lemma \ref{util}-(1)
implies
$c_1\subset \dom(\Pi)$.
Lemma \ref{util}-(2)
implies that $c_2=\Pi(c_1)$ is a $C^u$ curve
contained in $O'$. Definition \ref{sp}-(5)
implies
$L(c_2)\geq \epsilon^{-1}L(c_1)\geq R\cdot L(c_1)$.
Suppose we have constructed
a sequence $c_1,c_2,\cdots ,c_i$
of $C^u$ curves of
$O$ contained in $O'$ satisfying
$L(c_j)\geq R\cdot L(c_{j-1})$ for all
$2\leq j\leq i$.
If $c_i\cap l=\emptyset$ we
define $c_{i+1}=\Pi(c_i)$ and keep going.
If
$c_i\cap l\neq \emptyset$ we
define $c_{i+1}=Closure(\Pi(c_i^+))$.
In any case $c_{i+1}$ is a $C^u$ curve of $O$ contained
in $O'$.
In the first case we have
$L(c_{i+1})\geq \epsilon^{-1}L(c_i)\geq R\cdot L(c_i)$.
In the second case we have
$$
L(c_{i+1})=
L(\Pi(c_i^+))\geq \epsilon^{-1}L(c_i^+)
\geq (\epsilon^{-1}/2)\cdot L(c_i)
=R\cdot L(c_i).
$$
In this way we can construct an infinite
sequence $c_1,\cdots ,c_i,c_{i+1}, \cdots$
of $C^u$ curves of $O$ in $O'$ all of which satisfying
$L(c_{i+1})\geq R\cdot L(c_i)$.
It follows that
$$
L(c_i)\geq R^i\cdot L(c_1),
$$
for all $i$.
Since $l(c_1)>0$ and $R=(2\epsilon)^{-1}>1$
we conclude that
$$
\lim_{i\to\infty}L(c_i)=\infty.
$$
On the other hand,
let $S$ be the singular cross-section
containing $O$ given by Definition \ref{sp}-(1).
Let ${\cal F}^S$ be the projection of the
stable manifold in $U(\Lambda)$ over $S$.
As noted in Remark \ref{rk'} the leave
space $I^S$ of ${\cal F}^S$ is a finite union of
compact intervals. In particular
$I^S$ has finite diameter.
Since $O'\subset O\subset Int(S)$ we have that
all the curves $c_i$ are contained in $S$.
Since $c_i$ is a $C^u$ curve
we have by
Definition \ref{sp}-(6) that
$c_i$ have positive angle with the leaves of ${\cal F}^S$ (note that these leaves
are tangent to $E^s$).
So, we can project $c_i$ to
obtain an infinite sequence of intervals in $I^S$.
The lenght of these intervals goes to $\infty$
(as $i\to \infty$) since
$L(c_i)\to \infty$ (as $i\to \infty$).
This is a contradiction since $I^S$ has finite diameter.
This contradiction proves the result.
\end{proof}

{\flushleft{\bf Proof of Theorem \ref{thA}: }}
Let $\Lambda$ be a singular-hyperbolic
attractor on a compact $3$-manifold.
If $\Lambda$ has topological dimension
$<2$ then
$\Lambda$ would be one-dimensional
because it has regular orbits \cite{HW}.
It would follow from Proposition
\ref{p1} that $\Lambda$ has singular partitions with arbitrarily small size contradicting Theorem \ref{p2}. The proof follows.
\qed 

\section{Appendix}

In this section we state (and prove)
two technical lemmas which were
used in the proof of Theorem \ref{thA}.
These lemmas were proved in \cite{MP}
and here we reproduce these proofs
for the sake of completeness.
Let us state some definitions and notations.

First we define the Linear Poincar\'e Flow \cite{D}.
Let $X$ be a flow on a compact $3$-manifold $M$.
The Riemmanian Metric of $M$ is denoted
by $<\cdot,\cdot >$.
If $x$ is a regular point of $X$ (i.e. $X(x)\neq 0$),
we denote by $N_x=\{v_x\in T_xM: <v_x,X(x)>=0\}$ the orthogonal
complement of $X(x)$ in $T_xM$.
Denote $O_x:T_xM\to N_x$ the orthogonal projection onto $N_x$.
For every $t\in I\!\! R$ we define $P_x^t:N_x\to N_{X_t(x)}$
by
$$
P_x^t=O_{X_t(x)}\circ DX_t(x).
$$
It follows that $P=\{P_x^t:t\in I\!\! R,X(x)\neq 0\}$
satisfies the cocycle relation
$$
P^{s+t}_x=P^s_{X_t(x)}\circ P^s_x,
$$
for every $t,s\in I\!\! R$.
The parametrized family $P$ is called the {\em Linear Poincar\'e Flow} 
of $X$.

We denote by $\vol (v_x,w_x)$ the area of the parallelogram in $T_xM$ 
generated by $v_x,w_x\in T_xM$.
As $M$ is a compact manifold,
there is a constant $V\geq 1$ such that $V^{-1}\leq \vol (v_x,w_x)\leq V$,
$\forall x\in M$, $\forall v_x,w_x\in T_xM$ satisfying
$\| v_x\|=\| w_x\|=1$ and $<v_x,w_x>=0$.
For simplicity we shall assume that $V=1$.
In other words,
$$
\vol (v_x,w_x)=\| v_x\| \cdot \| w_x\|,
$$
$\forall x\in M$, $\forall v_x,w_x\in T_xM$ with $<v_x,w_x>=0$.

In addition,
$$
\vol(v_x,X(x))=\| O_xv_x\|\cdot\| X(x)\|,
$$
$\forall x\in M$ regular, $\forall v_x\in T_xM$.
In particular,
\begin{equation}
\label{e10}
\vol (DX_t(x)v_x,X(X_t(x)))=\| P_x^t(v_x)\|\cdot\| X(X_t(x))\|,
\end{equation}
$\forall x\in M$ regular, $\forall t\in I\!\! R$, $\forall v_x\in N_x$.

Recall that if $\Lambda$ is a singular hyperbolic set of $X$ with 
$(K,\lambda )$-splitting $T_\Lambda M=E^s_\Lambda\oplus E^c_\Lambda$,
then
$$
\mid
Det(DX_t/E^c_x)
\mid \geq K e^{\lambda t},
$$
$\forall x\in \Lambda$, $\forall t\geq 0$, where
$Det(\cdot)$ denotes the jacobian.
So, if $\Lambda$ is a singular hyperbolic set as above one has
$$
\vol (DX_t(x)v_x^c,DX_t(x)w_x^c)\geq K e^{\lambda t}\vol (v_x^c,w_x^c),
$$
$\forall x\in \Lambda$, $\forall t\geq 0$, $\forall v_x^c,w_x^c\in E^c_x$.

Remember that $U(\Lambda)$ denotes
a neighborhood of $\Lambda$ where
the splitting $E^s_\Lambda\oplus E^c_\Lambda$
extends to $E^s_{U(\Lambda)}\oplus
E^c_{U(\Lambda)}$.

\begin{lemma}
\label{AP1}
Let $\Lambda$ be a singular-hyperbolic
attractor of a $C^1$ flow $X$
on a compact $3$-manifold $M$.
Then for every $\alpha\in (0,1]$ there are a neighborhood
$U_{\alpha}\subset U(\Lambda)$ of $\Lambda$ and constants
$T_{\alpha},K_\alpha, \lambda_\alpha>0$ such that :
\begin{enumerate}
\item
If $x\in U_{\alpha}$ and $t\geq T_\alpha$, then
$$
DX_t(x)(C_\alpha(E^c_x))\subset C_{\frac{\alpha}{2}}(E_{X_t(x)}^c).
$$
\item
If $x\in U_{\alpha}$ is regular, $X(x)\in C_\alpha(E^c_x)$,
$t\geq T_\alpha$ and $v_x\in C_\alpha(E_x)\cap N_x$, then
$$
\| P^t_x(v_x)\| \cdot \|X(X_t(x))\| \geq K_\alpha e^{\lambda_\alpha t}
\cdot \| v_x\|\cdot \| X(x)\|.
$$
\end{enumerate}
\end{lemma}

\begin{proof}
Let $\Lambda$ and $\alpha\in (0,1]$ be as in the statement.
As mentioned above $T_{U(\Lambda)}M=E^s_{U(\Lambda)}
\oplus E^c_{U(\Lambda)}$ denotes the extension of
the $(K,\lambda)$-splitting
$T_\Lambda M=E^s_{\Lambda}\oplus E^c_{\Lambda}$ of $\Lambda$
to a neighborhood $U(\Lambda)$ of $\Lambda$.
Let $\pi^s$ the projection of $T_\Lambda$ on $E^s_\Lambda$,
and $\pi^c$ be the projection of $T_\Lambda$ on $E^c_\Lambda$.
Denote
$v_x=v_x^s+v^c_x\in E^s_x\oplus E^c_x=T_xM$
$\forall x\in U(\Lambda)$, $\forall v_x\in T_xM$.
In other words,
$v_x^s=\pi^s(v_x)$ and $v_x^c=\pi^c(v_x)$

As $E^s_\Lambda$ $(K,\lambda)$-dominates $E^c_\Lambda$ we have that
\begin{equation}
\label{nova}
\| DX_t(x)/E^s_x\| \leq K^{-1}e^{-\lambda t}m(DX_t(x)/E^c_x),
\end{equation}
$\forall x\in \Lambda$, $\forall t\geq 0$.

Fix $R>4$ such that
\begin{equation}
\label{e11}
\frac{K}{R}<1.
\end{equation}
Choosing $T^1=T^1_{\alpha}>0$ large enough one has
\begin{equation}
\label{et1}
\| DX_{T^1}(x)/E^s_x\|\leq\frac{K\alpha}{2R}m(DX_t(x)/E^c_x),
\end{equation}
$\forall x\in \Lambda$, $\forall t\geq 0$.

Since $E^s_\Lambda\oplus E^c_\Lambda$ is invariant we have
$\pi^{s(c)}_{X_t(x)}\circ DX_t(x)=DX_t(x)\circ \pi^{s(c)}_x$, and so
\begin{equation}
\label{nnova}
\angle(DX_t(x)v_x,E^c_{X_t(x)})
=
\frac{\| DX_t(x)v_x^s\|}{\| DX_t(x)v^c_x\|},\quad \forall x\in \Lambda
\quad \forall t\geq 0.
\end{equation}

Recall that $\angle$ denotes the tangent of the angle.
The inequality (\ref{et1}) and the last equality imply
$$
\angle(DX_{T^1}(x)v_x,E^c_{X_{T^1}(x)})\leq \frac{K\alpha}{2R},
\quad \forall x\in \Lambda, \,\,\forall v_x\in C_\alpha(E^c_x).
$$
So,
$$
DX_{T^1}(x)(C_\alpha(E^c_x))\subset
C_{\frac{K\alpha}{2R}}(E^c_{X_{T^1}}(x)),
\quad \forall x\in \Lambda.
$$

Choose a neighborhood $U^1=U_{\alpha}^1\subset U(\Lambda)$ of $\Lambda$ 
sufficiently close to $\Lambda$ such that
\begin{equation}
\label{e12}
DX_{T^1}(x)(C_\alpha(E^c_x))
\subset
C_{\frac{K\alpha}{R}}(E^c_{X_{T^1}}(x)),\quad \forall x\in U^1.
\end{equation}

On the other hand, using (\ref{nova}) we get 
$$
\frac{\|DX_t(x) v_x^s\|}{\|v_x^s\|}\leq
K^{-1}e^{-\lambda t}\frac{\|DX_t(x)v_x^c\|}{\|v_x^c\|}
$$
and so, 
$$\frac{\|DX_t(x) v_x^s\|}{\|DX_t(x)v_x^c\|}\leq
K^{-1}e^{-\lambda t}\frac{\|v_x^s\|}{\|v_x^c\|}=
K^{-1}e^{-\lambda t}\angle (v_x,E^c_x).
$$
So, by (\ref{nnova}), we get
$$
\angle(DX_r(x)v_x,E^c_{X_r(x)})\leq
%\frac{\| DX_r(x)v_x^s\|}{
%\| DX_r(x)v_x^c\|}
%\leq
K^{-1}e^{-\lambda r}\angle(v_x,E^c_x)
\leq K^{-1}\angle(v_x,E^c_x),
$$
$\forall x\in \Lambda$, $\forall r\in [0,T^1]$, $\forall v_x\in T_xM$.
This implies
$$
DX_r(x)\left(C_{\frac{K\alpha}{R}}(E^c_{x})\right)
\subset C_{\frac{\alpha}{R}}(E^c_{X_r(x)}),
\quad \forall x\in \Lambda \quad\forall r\in [0,T^1].
$$
Choose a neighborhood $V^2=V^2_{\alpha}\subset U^1$
of $\Lambda$ sufficiently close to $\Lambda$ such that
\begin{equation}
\label{et1.5}
DX_r(x)\left(C_{\frac{K\alpha}{R}}(E^c_{x})\right)\subset
C_{\frac{2\alpha}{R}}(E^c_{X_r(x)}),\,\,
\forall x\in V^2, \,\,\forall r\in [0,T^1].
\end{equation}
As $\Lambda$ is an attractor there is a neighborhood $U^2\subset V^2$
of $\Lambda$ such that
$$
X_t(U^2)\subset V^2, \quad \forall t\geq 0.
$$
Now, let $x\in U^2$ and $t\geq T^1$ be given.
Then, $t=nT^1+r$ for some integer $n\geq 1$ and some $r\in [0,T^1]$.
Thus,
$$
DX_t(x)(C_\alpha(E^c_x))=DX_{nT^1+r}(x)(C_\alpha(E^c_x))=
$$
\begin{equation}
\label{star}
=DX_r(X_{nT^1}(x))(DX_{nT^1}(x)(C_\alpha(E^c_x))).
\end{equation}
Using  (\ref{e12}) and (\ref{star}) recursively, and that  
$\frac{K\alpha}{R}<\alpha$, and $n\geq 1$, we obtain
$$
DX_{nT^1}(x)(C_\alpha(E^c_x))=DX_{(n-1)T^1}(X_{T^1}(x))
(DX_{T^1}(x)C_\alpha(E^c_x))\subset
$$
$$
\subset
DX_{(n-1)T^1}(X_{T^1}(x))\left(C_{\frac{K\alpha}{R}}(E^c_{X_{T^1}(x)})
\right)\subset DX_{(n-1)T^1}(X_{T^1}(x))(C_\alpha(E^c_{X_{T^1}(x)}))
\subset
$$
$$
\subset
\cdots\subset DX_{T^1}(X_{(n-1)T^1}(x))(C_\alpha(E^c_{X_{(n-1)T^1}(x)}))
\subset C_{\frac{K\alpha}{R}}(E^c_{X_{nT^1}(x)}).
$$
%because $\frac{K\alpha}{R}<\alpha$ and $n\geq 1$.
Henceforth
$$
DX_{nT^1}(x)(C_\alpha(E^c_x))\subset C_{\frac{K\alpha}{R}}(E^c_{X_{nT^1}(x)}).
$$
Applying $DX_r(X_{nT^1}(x))$ to both sides of the last expression, 
replacing in (\ref{star}) and using (\ref{et1.5})  we obtain
$$
DX_t(x)(C_\alpha(E^c_x))\subset DX_r(X_{nT^1}(x))\left(
C_{\frac{K\alpha}{R}}(E^c_{X_{nT^1}(x)})\right)\subset
C_{\frac{2\alpha}{R}}(E^c_{X_t(x)}).
$$
As $R>4$ we have $\frac{2\alpha}{R}< \frac{\alpha}{2}$ and so
\begin{equation*}
\label{et2}
DX_t(x)(C_\alpha(E^c_x))
\subset
C_{\frac{\alpha}{2}}(E^c_{X_t(x)}).
\end{equation*}
$\forall x\in U^2$, $\forall t\geq T^1$,
proving (1) of Lemma \ref{AP1}.

Throughout we fix the neighborhood $U^2$ of $\Lambda$ and the constant
$T^1>0$ obtained above.

As $E^c_\Lambda$ is $(K,\lambda)$-volume expanding we have
$$
\vol (DX_t(x)v_x^c,DX_t(x)w_x^c) \geq Ke^{\lambda t}\vol (v_x^c,w_x^c),
$$
$\forall x\in \Lambda$, $\forall t\geq 0$, $\forall v_x^c,w_x^c\in  E^c_x$.

Clearly there is $L>1$ such that
$$
L^{-1}\cdot \vol (v_x^c,w_x^c)\leq\vol (v_x,w_x)\leq L\cdot \vol (v_x^c,w_x^c),
$$
$\forall x\in \Lambda$, $\forall v_x,w_x\in C_\alpha(E^c_x)$,
$\forall \alpha\in (0,1]$.
Applying the last two relations and the invariance of
$E^s_\Lambda\oplus E^c_\Lambda$ we obtain
$$
\vol (DX_t(x)v_x,DX_t(x)w_x)\geq L^{-1} \vol (DX_t(x)v^c_x,DX_t(x)w_x^c)\geq
$$
$$
\geq L^{-1}Ke^{\lambda t} \vol (v_x^c,w_x^c)
\geq L^{-2}Ke^{\lambda t}\vol (v_x,w_x),
$$
$\forall x\in \Lambda$, $\forall t\geq T^1$, 
$\forall v_x,w_x\in C_\alpha(E^c_x)$
(note that $DX_t(x)v_x, DX_t(x)w_x\in C_\alpha(E^c_{X_t(x)})$ since 
$t\geq T^1$).

Choose $S>0$ large so that
\begin{equation}
\label{e13}
\frac{S}{L^{-2}K}>1.
\end{equation}
It follows that there is $T^2=T^2_{\alpha}>T^1$ such that
$$
\vol (DX_{T^2}(x)v_x,DX_{T^2}(x)w_x)
\geq
\frac{2S}{L^{-2}K}\vol(v_x,w_x),
$$
$\forall x\in \Lambda$,
$\forall v_x, w_x\in C_\alpha(E^c_x)$.
In particular,
$$
\inf\{\vol(DX_{T^2}(x)v_x,DX_{T^2}(x)w_x): x\in \Lambda,
$$
$$
v_x,w_x\in C_\alpha(E^c_x),
\| v_x\|=\| w_x\|=1, <v_x,w_x>=0\}\geq \frac{2S}{L^{-2}K}.
$$
 
Since $\Lambda$ is compact there is a neighborhood 
$V^3=V^3_{\alpha}\subset U^2$ of $\Lambda$ so that
$$
\inf\{
\vol(DX_{T^2}(x)v_x,DX_{T^2}(x)w_x): x\in U^3,
$$
$$
v_x,w_x\in C_\alpha(E^c_x),\| v_x\|=\| w_x\|=1,
<v_x,w_x>=0\}\geq \frac{S}{L^{-2}K}.
$$
Then,
\begin{equation}
\label{e40}
\vol(DX_{T^2}(x)v_x,DX_{T^2}(x)w_x)
\geq \frac{S}{L^{-2}K} \vol (v_x,w_x),
\end{equation}
$\forall x\in U^3$, $\forall v_x,w_x\in C_\alpha(E^c_x)$ with $<v_x,w_x>=0$.
As $\Lambda$ is an attractor there is a neighborhood $U^3\subset V^3$
of $\Lambda$ such that
$$
X_t(U^3)\subset V^3,
\quad
\forall t\geq 0.
$$
We have
$$
\|P^{nT^2}_xv_x\| \,\,\| X(X_{nT^2}(x))\|=
\|P^{T^2}_{X_{(n-1)T^2}(x)}(P^{(n-1)T^2}_x v_x\|\,\,
\|X(X_{T^2}(X_{(n-1)T^2}(x))\|.
$$
Call $z=X_{(n-1)T^2}(x)$, and $v_z=P_x^{(n-1)T^2}v_x$.
From the last equality, using that 
$X(X_{nT^2}(x))=DX_{nT^2}(x)(X(x))$,  
$v_z$ is orthogonal to $z$,
and combining (\ref{e10}) with (\ref{e40}) we get
$$
\|P^{nT^2}_xv_x\| \,\,\| X(X_{nT^2}(x))\|=
\|P_z^{T^2}v_z\|\,\,\|X(X_{T^2}(z))\|=
$$
$$
\vol(DX_{T^2}(z)v_z, X(X_{T^2}(z)))=
\vol(DX_{T^2}(z)v_z, DX_{T^2}(z)(X(z))
\geq
$$
$$
\frac{S}{L^{-2}K}\vol(v_z,X(z))=
\frac{S}{L^{-2}K}\vol(P^{(n-1)T^2}_xv_x, X(X_{(n-1)T^2}(x))).
$$
Thus,
\begin{equation}
\label{et3}
\|
P^{nT^2}_xv_x\|\cdot \| X(X_{nT^2}(x))\|\geq
\left\{
\frac{S}{L^{-2}K}
\right\}^n
\| v_x\|\cdot \| X(x)\|,
\end{equation}
$\forall x\in U^3$ regular with $X(x)\in C_\alpha(E^c_x)$,
$\forall n\in I\!\! N$, $\forall v_x\in C_\alpha(E^c_x)\cap N_x$
(recall that $N_x$ denotes the orthogonal complement of $X(x)$ in $T_xM$).

On the other hand
$$
\vol(DX_r(x)v_x,DX_r(x)w_x) \geq L^{-2}K\cdot \vol(v_x,w_x),
$$
$\forall x\in \Lambda,\quad
\forall r\in [0,T^2],
\quad\forall v_x,w_x\in C_\alpha(E^c_x).
$

As before there is a neighborhood $V^4=V^4_{\alpha}\subset U^3$ of
$\Lambda$ such that
\begin{equation}
\label{e41}
\vol(DX_r(x)v_x,DX_r(x)w_x) \geq \frac{L^{-2}K}{2}\vol(v_x,w_x),
\end{equation}
$\forall x\in V^4$,
$\forall v_x,w_x\in C_\alpha(E^c_x)$
with $<v_x,w_x>=0$, $r\in [0,T^2]$.
As $\Lambda$ is an attractor there is a neighborhood
$U^4\subset V^4$ of $\Lambda$ such that
$$
X_t(U_4)\subset V^4,
\quad \forall t\geq 0.
$$

Now, let
$x\in U^4$ regular with $X(x)\in C_\alpha(E^c_x)$, $t\geq T^2$ and
$v_x\in C_\alpha(E^c_x)\cap N_x$.
Then, $t=nT^2+r$ for some integer $n\geq 1$ and some $r\in [0,T^2]$.

Applying (\ref{et3}), (\ref{e41}), and using (\ref{e10}) and (\ref{e13}) we obtain
$$
\| P^t_xv_x\|\cdot\|X(X_t(x))\|=
\| P^r_{X_{nT^2}(x)}P^{nT^2}_xv_x\|
\cdot
\|DX_r(X_{nT^2}(x))X(X_{nT^2}(x))\|
=
$$
$$
=\vol\left(DX_r(X_{nT^2}(x))P^{nT^2}_xv_x,
DX_r(X_{nT^2}(x))X(X_{nT^2}(x))
\right)\geq
$$
$$
\geq \frac{L^{-2}K}{2}
\vol(P^{nT^2}_xv_x,X(X_{nT^2}(x)))=
\frac{L^{-2}K}{2}
\|
P^{nT^2}_xv_x\|\cdot
\|
X(X_{nT^2}(x))\|
\geq
$$
$$
\geq
\left(\frac{L^{-2}K}{2}\right)\cdot
\left\{
\frac{SN}{L^{-2}K}
\right\}^n\cdot
\| v_x\|\cdot \|
X(x)\|=
$$
$$
=
\left(\frac{L^{-2}K}{2}\right)\cdot
\left(
\frac{S}{L^{-2}K}
\right)^{-\frac{r}{T^2}}
\left\{
\left(\frac{S}{L^{-2}K}\right)^{\frac{1}{T^2}}
\right\}^t
\cdot\| v_x\|\cdot \|
X(x)\|
\geq
$$
$$
\geq
\left(\frac{L^{-4}K^2}{2S}\right)\cdot
\left\{
\left(\frac{S}{L^{-2}K}\right)^{\frac{1}{T^2}}
\right\}^t\cdot
\| v_x\|\cdot \|
X(x)\|.
$$
%The lemma follows from (\ref{et2})
%and the last computation by 
Thus, choosing $U_{\alpha}=U^4$, $T_{\alpha}=T^2$, 
$K_\alpha=\frac{L^{-4}K^2}{2S}$
and $\lambda_\alpha=\ln\left(\frac{S}{L^{-2}K}\right)^{\frac{1}{T^2}}>0$
we obtain (2) of Lemma \ref{AP1}.
\end{proof}

\begin{lemma}
\label{AP2}
Let $\Lambda$ a singular-hyperbolic attractor
of $X$, $\alpha>0$ and
$U_\alpha \subset U(\Lambda)$ be a neighborhood of $\Lambda$.
Let $T_\Lambda M=E^s_\Lambda
\oplus E^c_\Lambda$ the splitting of $\Lambda$ and $E^c_{U_\alpha}$
be a continuous extension of $E^c_\Lambda$ to $U_\alpha$.
Then, for every singularity $\sigma$ of $X$ in $\Lambda$
there are singular cross-sections $S^t,S^b$ associated to $\sigma$
such that $S^t\cup S^b\subset U_\alpha$,
$$
\Lambda\cap (\partial^hS^t\cup\partial^hS^b)=
\emptyset,\,\,\mbox{ and }\,\,
X(x)\in C_\alpha(E^c_x),
$$
$\forall x\in S^t\cup S^b$.
\end{lemma}
\begin{proof}
Let $\Lambda,\alpha$ be as in the statement.
By Remark \ref{rk} applied to $U=U_\alpha$
there are singular cross-sections
$S^t_0,S^b_0$ associated to $\sigma$ such that
$S^t_0,S^b_0\subset U(\Lambda)$ and
$\Lambda\cap (\partial^hS^t_0\cup\partial^hS^b_0)=\emptyset$. 
Recall that $U(\Lambda)$ is the neighborhood of $\Lambda$ where
the $(K,\lambda)$-splitting $T_\Lambda M=E^s_\Lambda \oplus E^c_\Lambda$
has an extension to $T_{U(\Lambda)}M=E^s_{U(\Lambda)}\oplus E^c_{U(\Lambda)}$,
with $E^s_{U(\Lambda)}$ invariant and contracting.
We denote by $l^t_0,l^b_0$ the singular curves of $S^t_0,S^b_0$
respectively.

Choose two sequences of
singular cross-sections
$S^t_n\subset S^t_0,S^b_n\subset S^b_0$ associated to $\sigma$
satisfying
\begin{description}
\item{(a)}
$\Lambda\cap (\partial^hS^t_n\cup\partial^hS^b_n)=
\emptyset$;
\item{(b)}
$\diam(S^t_n),\diam(S^b_n)\to 0$ as $n\to \infty$;
\item{(c)}
the singular curves
$l^t_n,l^b_n$ of $S^t_n,S^b_0$ satisfy
$l^t_n=l^t_0$, $l^b_n=l^b_0$,
$\forall n$.
\end{description}

The properties (b) and (c) imply that $\forall n$ large there is $T=T_n>0$ such that 
$$
\angle(X(X_T(x)), E^c_\sigma)<\frac{\alpha}{2},
$$
$\forall x\in S^t_n\cup S^b_n$.
As $E^c_{U_\alpha}$ is a continuous
extension of $E^c_\Lambda$, we have
that $E^c_{X(X_T(x))}$ is close to
$E^c_\sigma$, $\forall x\in S^t_n\cup S^b_n$.
Then,
$$
X_T(x)\in U_\alpha,\,\mbox{ and }\,\,
X(X_T(x))\in C_\alpha(E^c_{X(X_T(x))}),
$$
$\forall n$ large, $\forall x\in S^t_n\cup S^b_n$.

By the Property (a) and the above relation
we have that for $n$ large enough
$S^t=X_T(S^t_n)$ and $S^b=X_T(S^b_n)$
are singular cross-sections associated
to $\sigma$ satisfying the required properties.
\end{proof}

\vspace{0.2cm}
\noindent C. A. Morales\\
Instituto de Matem\'atica \\
Universidade Federal do Rio de Janeiro \\
C. P. 68.530, CEP 21.945-970 \\
Rio de Janeiro, R. J. , Brazil \\
e-mail: {\em morales@impa.br}

\end{document}